\documentclass[reqno,12pt]{amsart}

\usepackage{amsmath,amsthm,amsfonts,amssymb} \date{}
\newtheorem{theorem}{Theorem}[section]
\newtheorem{lemma}{Lemma}[section]

\newtheorem{remark}{Remark}[section]
\newtheorem{example}{Example}[section]
\numberwithin{equation}{section}

\begin{document}

\title[Global regularity and bounds]
{Global regularity and bounds for solutions of
 parabolic equations for probability measures}

\author[V.I.~Bogachev]{Vladimir~I.~Bogachev}
\address{Department of Mechanics and Mathematics,
         Moscow State University, 119992 Moscow, Russia}

\author[M.~R\"ockner]{Michael~R\"ockner}
\address{Department of Mathematics, Purdue
University, 150 N.~University Str.,
West Lafay\-ette, IN 47907-2067, USA}

\author[S.V.~Shaposhnikov]{Stanislav~V.~Shaposhnikov}
\address{Department of Mechanics and Mathematics,
         Moscow State University, 119992 Moscow, Russia}

\thanks{
This work has been supported by the projects RFBR
04-01-00748,
the Scientific Schools Grant 1758.2003.1,
 DFG 436 RUS 113/343/0(R),
 INTAS 03-51-5018,
and the  SFB 701 at the University of Bielefeld.}

\maketitle

{\small
Given a second order parabolic operator
$$
Lu(t,x):=\frac{\partial u(t,x)}{\partial t}+
a^{ij}(t,x)\partial_{x_i}\partial_{x_j}u(t,x)+b^i(t,x)\partial_{x_i}u(t,x),
$$
we consider the weak parabolic equation $L^{*}\mu=0$
for Borel probability measures on $(0,1)\times\mathbb{R}^d$.
The equation is understood as the equality
$$
\int_{(0,1)\times\mathbb{R}^d} Lu\, d\mu =0
$$
for all smooth functions $u$ with compact support
in~$(0,1)\times\mathbb{R}^d$.
This equation is satisfied for the transition
probabilities of the diffusion process associated
with~$L$.
We show that under broad assumptions $\mu$ has the form
$\mu=\varrho(t,x)\, dt\, dx$, where the function
$x\mapsto \varrho(t,x)$ is  Sobolev,
$|\nabla_x \varrho(x,t)|^2/\varrho(t,x)$ is Lebesgue
integrable over $[0,\tau]\times\mathbb{R}^d$, and
$\varrho\in L^p([0,\tau]\times\mathbb{R}^d)$
for all $p\in [1,+\infty)$ and $\tau<1$.
Moreover, a sufficient condition for the uniform
boundedness of $\varrho$ on $[0,\tau]\times\mathbb{R}^d$
is given.


}

\section{Introduction and notation}

The objective of this work is
to give efficient conditions for the global
Sobolev regularity and integrability of densities of solutions
of the
parabolic equations of the form
\begin{equation}\label{e1.1}
L^{*}\mu =0
\end{equation}
for Borel measures $\mu$ on $(0,1)\times\mathbb{R}^d$.
Such equations have been recently investigated in
 \cite{BDPR}, \cite{BDPR05}, \cite{BKR01}, \cite{Stannat}.
Here $L$ is a second order parabolic operator
$$
Lu(t,x):=\frac{\partial u(t,x)}{\partial t}+
a^{ij}(t,x)\partial_{x_i}\partial_{x_j}u(t,x)+b^i(t,x)\partial_{x_i}u(t,x),
$$
and the interpretation of our equation is the following.
We shall say that a Borel probability
measure $\mu$ on $(0,1)\times\mathbb{R}^d$
represented in the form $\mu(dt\, dx)=\mu_t(dx)\, dt$ by means
of a family of Borel measures $(\mu_t)_{t\in [0,1)}$
on $\mathbb{R}^d$
satisfies the weak parabolic equation (\ref{e1.1})
if the functions $a^{ij}$ and $b^i$ are
integrable on every compact set in $(0,1)\times \mathbb{R}^d$
with respect to the measure $\mu=\mu_t\, dt$ and,
for every $u\in C_0^\infty ((0,1)\times \mathbb{R}^d)$, one has
\begin{equation}\label{e1.2}
\int_{(0,1)\times\mathbb{R}^d} Lu\, d\mu=
\int_0^1 \int_{\mathbb{R}^d} Lu(t,x)\mu_t(dx)\, dt =0 .
\end{equation}
We shall say that $\mu$ satisfies the initial condition $\mu_0=\nu$ at $t=0$
if $\nu$ is a Borel measure on $\mathbb{R}^d$ and
\begin{equation}\label{e1.3}
\lim\limits_{t\to 0}\int_{\mathbb{R}^d} \zeta (x)\, \mu_t(dx)=
\int_{\mathbb{R}^d} \zeta (x)\, \nu(dx)
\end{equation}
for all $\zeta\in C_0^\infty (\mathbb{R}^d)$.

The same definitions are introduced in the case where
$\mathbb{R}^d$ is replaced by an open set $\Omega\subset\mathbb{R}^d$
or by an open set in a Riemannian manifold.

Equation (\ref{e1.1}) is satisfied for the transition
probabilities of the diffusion process governed
by the stochastic differential equation
$$
d\xi_t=\sqrt{2A(t,\xi_t)} dw_t+ b(t,\xi_t)dt
$$
provided that such a diffusion
exists and the coefficients $A$ and $b$ satisfy certain conditions.
However, (\ref{e1.1}) can be considered regardless of any
probabilistic assumptions. Moreover, a study of this equation
in a purely analytical setting may be useful for constructing
an associated diffusion (see \cite{Stannat}).

Our main result states that the density $\varrho$ of any
solution has the property that
$\varrho(t,\,\cdot\,)$ is Sobolev on $\mathbb{R}^d$ and
$|\nabla\varrho(t,x)|^2/\varrho(t,x)$ is integrable over
$[0,\tau]\times\mathbb{R}^d$ provided that the functions
$|b|$ and $\ln(|x|+1)$ are in~$L^2(\mu)$, the coefficient
$A$ is uniformly bounded, uniformly invertible and uniformly
Lipschitzian in~$x$, and the initial distribution
$\mu_0=\varrho(0,\,\cdot\,)\, dx$ has
finite entropy. The assumptions on $A$ can be relaxed if
$b$ has certain additional local integrability.
An efficient condition in terms of Lyapunov functions is given
in order to ensure the square integrability of $|b|$
and $\ln(|x|+1)$ with respect to the solution~$\mu$.
The main result enables us to show that $\varrho$ belongs to all
$L^p([0,\tau]\times \mathbb{R}^d)$ whenever
$\tau<1$, provided that
 $\sup_t \|b(t,\,\cdot\,)\|_{L^d(\mu_t)}<\infty$
and $\varrho(0,\,\cdot\,)\in L^p(\mathbb{R}^d)$
for all~$p\ge 1$. If
$|b|\in L^\beta(\mu)$
for some $\beta>d+2$ and
$\varrho(0,\,\cdot\,)\in L^\infty(\mathbb{R}^d)$, then
the density $\varrho$ is uniformly bounded on
$[0,\tau]\times \mathbb{R}^d$ whenever $\tau<1$.
 By using this assertion we obtain pointwise upper bounds
 of the form $\varrho(t,x)\le \Phi(x)^{-1}$.
Note that unlike many known results on the global boundedness
of solutions, it is not required here that the drift term
be dissipative or potential.

Analogous results in the elliptic case have been obtained
in \cite{BR95}, \cite{BKR96}, \cite{BRW},
\cite{MPR}, \cite{BKR05a} and~\cite{BKR05}.
One might regard the elliptic case as the situation
when the solution and the coefficients are independent of time.
Then our parabolic result does not recover the elliptic one, because
the initial distribution (which in this case coincides with the solution)
must have finite entropy, and the latter assumption
cannot be completely removed in the parabolic case.
On the other hand,  a reasonable parabolic analogue of the
elliptic result  might be as follows: the integrability of
$|\varrho|^2/\varrho$ on $[\tau_1,\tau_2]\times\mathbb{R}^d$ for any
closed interval $[\tau_1,\tau_2]\subset (0,1)$ without restrictions
on the initial distribution. So far
we have not succeeded in investigating this second possibility.
Of course, if for some $\tau>0$ the measure $\mu_\tau$ has finite
entropy, then our hypotheses are satisfied on~$[\tau,1]$.

Our result will be applied in a forthcoming paper on the
uniqueness problem for parabolic equations for measures.
It can also be  useful in the study of transition probabilities
of diffusion processes and in Nelson's dynamics
(see \cite{Carlen}, \cite{CatF},
\cite{CattiauxL}, \cite{Foellmer}, \cite{Nelson}).

Let $W^{p,1}(\mathbb{R}^d)$ denote the Sobolev space
of functions that belong to $L^p(\mathbb{R}^d)$ with
their generalized partial derivatives.
This space is equipped  with the standard norm
$$
\|f\|_{W^{p,1}}:=\|f\|_p+\|\nabla f\|_p,
$$
where $\|\,\cdot\,\|_p$ denotes the $L^p(\mathbb{R}^d)$-norm
on scalar or vector functions.
The symbol $L^{p,q}$, where $1\le p,q<\infty$,
will stand for the
space of all
measurable functions $f$ on $[0,1]\times\mathbb{R}^d$ with
finite norm
$$
\|f\|_{p,q}:=\Bigl(\int_0^1\Bigl(\int_{\mathbb{R}^d}
|f(t,x)|^p\, dx\Bigr)^{q/p}\, dt\Bigr)^{1/q}.
$$
The space $L^{p,\infty}$, where $1\le p<\infty$,
consists of all measurable functions
$f$ on $[0,1]\times\mathbb{R}^d$ with
$t\mapsto \|f(t,\,\cdot\,)\|_{L^p(\mathbb{R}^d)}\in
L^\infty [0,1]$.
For analogous spaces of functions on
$[0,\tau]\times\mathbb{R}^d$ we used the notation
 $L^{p,q}([0,\tau]\times\mathbb{R}^d)$.
Finally, $\mathbb{H}^{p,q}([0,\tau]\times\mathbb{R}^d)$ denotes the
space of all
measurable functions $f$ on $[0,\tau]\times\mathbb{R}^d$ with
finite norm
$$
\|f\|_{\mathbb{H}^{p,q}([0,\tau]\times\mathbb{R}^d)}:
=\Bigl(\int_0^\tau \|f(t,\,\cdot\,)\|_{W^{p,1}}^q\,
dt\Bigr)^{1/q}.
$$

For simplicity of notation
the gradient of a function $u$ on $(0,1)\times\mathbb{R}^d$
with respect to the argument from $\mathbb{R}^d$ is
 denoted by~$\nabla u$, i.e.,
 $$
 \nabla u(t,x):=\nabla _x u(t,x)=
 (\partial_{x_i}u(t,x),\ldots,\partial_{x_d}u(t,x)).
 $$

We use the standard rule of summation with respect
to repeated indices, e.g.,
$$
\partial_{x_i} a^{ij}:=\sum_{i=1}^d \partial_{x_i} a^{ij},
\quad
a^{ij}\partial_{x_i}\partial_{x_j}u
:=\sum_{i,j=1}^d  a^{ij}\partial_{x_i}\partial_{x_j}u.
$$

We say that a nonnegative measure $\mu_0$ on $\mathbb{R}^d$ has finite
entropy if $\mu_0=\varrho_0\, dx$ and
$\varrho_0\ln\varrho_0\in L^1(\mathbb{R}^d)$,
where we set $0\ln 0:=0$.
The entropy of $\varrho_0$ is the integral of~$\varrho_0\ln\varrho_0$.

It was shown in \cite{BKR01} that if the coefficient
$A$ is nondegenerate, then
$\mu$ is absolutely continuous with respect to Lebesgue measure
on $(0,1)\times\mathbb{R}^d$. The corresponding density
will be denoted by~$\varrho$.

A sufficient condition for the existence of a solution
in the class of probability measures is the following
(see \cite{BDPR05} which improves~\cite{BDPR}).
The coefficients $a^{ij}$ and $b^i$ are Borel functions on
$[0,1]\times\mathbb{R}^d$ such that $A(t,x)=(a^{ij}(t,x))$ is a
nonnegative symmetric matrix and there is $p>d+2$ such that, for every
ball $B$ and all $i,j\le d$, one has

{\rm(C1)} $\inf_{(t,x)\in [0,1]\times B}
\det A(t,x)\ge M_1(B)>0$ and
$\sup _{t\in [0,1]} \|a^{ij}(t,\,\cdot\,)\|_{W^{p,1}(B)}\le M_2(B)$,

{\rm(C2)} $\sup _{t\in [0,1]} \|b^{i}(t,\,\cdot\,)\|_{L^{p}(B)}\le M_3(B)$.

If (C1) and (C2) are fulfilled and there is a nonnegative
function $V$ on $\mathbb{R}^d$ such that
$\lim\limits_{|x|\to+\infty} V(x)=+\infty$ and
$LV\le c_1V+c_2$ for some constants $c_1$ and $c_2$, then
a solution $\mu=\mu_t(dx)\, dt$ with probability measures $\mu_t$
exists for every initial distribution~$\mu_0$
and the function
$$
 t\mapsto \int_{\mathbb{R}^d} \zeta(x)\, \mu_t(dx)
 $$
 is continuous on $[0,1)$ for every~$\zeta\in C_0^\infty(\mathbb{R}^d)$.
Moreover, if $V\in L^1(\mu_0)$, then
$$
\int_0^1\int_{\mathbb{R}^d} V(x)\, \mu_t(dx)\, dt<\infty.
$$
In addition, if $LV\le c_2$, then
$$
\int_{\mathbb{R}^d} V(x)\mu_t(dx)\le c_2
$$
for almost every $t$.
For example, if the coefficient $A$ is uniformly bounded
and
$$
\langle b(t,x),x\rangle \le k_1|x|^2+k_2,
$$
then we take the function $V(x)=\ln (|x|^2+1)$.
This gives an estimate
$LV\le const$, hence the integrals of $\ln(|x|+1)$
against $\mu_t$ are uniformly bounded
provided that $\ln(|x|+1)\in L^1(\mu_0)$.
If the coefficient $A$ is uniformly bounded
and
$$
\langle b(t,x),x\rangle \le k_1|x|^2\ln(|x|+1)+k_2,
$$
then we set $V(x)=|\ln (|x|^2+1)|^2$ and obtain
$LV\le c_1V+c_2$, hence $|\ln (|x|+1)|^2$ is $\mu$-integrable
provided that it is $\mu_0$-integrable.

Our main estimate will be established in the two cases corresponding
to two different approaches:

1)~when (C1) is replaced by a stronger assumption and (C2) is
replaced by the condition that
$|b|,\, \ln \max(|x|,1)\in L^2(\mu)$,

2)~the condition
$|b|,\, \ln \max(|x|,1)\in L^2(\mu)$
is imposed in addition to (C1), (C2) and
a certain global condition on~$A$.

Set
$$
\Theta_A(t,x):=\sum\limits_{j=1}^d \Bigl|
\sum_{i=1}^{d}\partial_{x_i}a^{ij}(t,x)\Bigr| .
$$

\section{Bounds on logarithmic gradients}

Our first main result establishes the square integrability of
the  logarithmic gradient of~$\mu$, i.e., the mapping
$\nabla\varrho/\varrho$, with respect to~$\mu$.
If $\varrho(t,\,\cdot\,)\in W^{1,1}_{loc}$, then we use the following convention:
$\nabla\varrho(t,x)/\varrho(t,x):=0$ if $\varrho(t,x)=0$.

\begin{theorem}\label{t2.1}
Suppose $\mu$, where each $\mu_t$ is a probability
measure, satisfies {\rm(\ref{e1.1})},~{\rm(\ref{e1.3})}.
 Let

{\rm(i)} the mapping
$A$ be uniformly bounded with
$A(t,x)\ge \alpha\cdot I$ for some constant~$\alpha>0$,
and let the functions $x\mapsto a^{ij}(t,x)$ be Lipschitzian
with constant~$\lambda$,

 {\rm(ii)} $|b|\in L^2(\mu)$.

\noindent
Assume also
that the function $\Lambda(x):=\ln \max(|x|,1)$ is in $L^2(\mu)$
{\rm(}which is the case if, e.g., $\langle b(t,x),x\rangle
\le C_1|x|^2\Lambda(x)+C_2$ with some constants $C_1$
and~$C_2$ and $\Lambda\in L^2(\mu_0)${\rm)}.
If $\mu_0$ has finite entropy, then
$\mu_t=\varrho(t,\,\cdot\,)\, dx$, where $\varrho(t,\,\cdot\,)\in
W^{1,1}_{loc}$, and for each $\tau<1$ one has
\begin{equation}\label{e2.1}
\int_0^\tau \int_{\mathbb{R}^d}
\frac{|\nabla \varrho(t,x)|^2}{\varrho(t,x)}\, dx\, dt <\infty .
\end{equation}
In particular, we have
$\sqrt{\varrho}\in \mathbb{H}^{2,2}
([0,\tau]\times\mathbb{R}^d))$ and
$\varrho\in L^{d/(d-2),1}([0,\tau]\times\mathbb{R}^d))$ if $d>2$, and
$\varrho\in L^{s,1}([0,\tau]\times\mathbb{R}^d))$ for all $s\in [1,\infty)$ if~$d=2$.

If the integrals ${\displaystyle \int_{\mathbb{R}^d} \varrho(t,x)
\Lambda(x)\, dx}$ are bounded as $t\to 1$
 {\rm(}which is the case, e.g.,
 if $\langle b(t,x),x\rangle\le C_1|x|^2+C_2$
 with some constants $C_1$ and~$C_2$ and $\Lambda\in L^1(\mu_0)${\rm)}, then
{\rm(\ref{e2.1})} is true for~$\tau=1$.
\end{theorem}
\begin{proof}
We shall use the following fact (see, e.g., \cite[Lemma~2.1]{BKR96}):
 given two nonnegative functions $f_1,f_2\in L^1(\mathbb{R}^d)$,
 for any measurable function
$\psi$ such that $|\psi|^2f_1\in L^1(\mathbb{R}^d)$
one has
\begin{equation}\label{e2.2}
\int_{\mathbb{R}^d}
\frac{|(\psi f_1)*f_2|^2}{f_1*f_2}\, dx
\le
\int_{\mathbb{R}^d}|\psi|^2f_1\, dx
\int_{\mathbb{R}^d}f_2\, dx,
\end{equation}
where $|(\psi f_1)*f_2(x)|^2/(f_1*f_2(x)):=0$ if $f_1*f_2(x)=0$.

For a function $w\in C_0^\infty(\mathbb{R}^d)$ we set
$$
\varrho *w(t,x):=\int_{\mathbb{R}^d}
w(x-y)\varrho(t,y)\, dy,\quad x\in\mathbb{R}^d .
$$
Here and
in what follows the convolutions are always taken
with respect to the variable from~$\mathbb{R}^d$.
Equation  (\ref{e1.2}) and the inclusion $|b|\in L^2(\mu)$
yield that the following
equality holds in Sobolev's sense:
\begin{equation}\label{e2.3}
\partial_t(\varrho * w)=
(a^{ij}\varrho)
*\partial_{x_i}\partial_{x_j}w-(b^i\varrho)*\partial_{x_i}w.
\end{equation}
We shall deal with a version of $\varrho * w$
(denoted by the same symbol) defined by the formula
  \begin{equation}\label{e2.4}
 \varrho *w(t,x):=\varrho*w(0,x)+\int_0^t v(s,x)\, ds ,
 \end{equation}
where $v$ is the right-hand side of~(\ref{e2.4}).
Since $|b|\in L^2(\mu)$ and the functions $a^{ij}$ are bounded,
one has $v\in L^1([0,1]\times\mathbb{R}^d)$.
Hence the function $\varrho *w$
is absolutely continuous in $t$ on $[0,1]$ and belongs
to the class $C_b^\infty(\mathbb{R}^d)$ in~$x$.
For almost every~$t$, the indicated version coincides for all
$x$ with the initial version defined by the convolution.
It will be important as well that this is true for~$t=0$.
Since the initial version does not exceed $\sup_x |w(x)|$ in the
absolute value, the same is true for our new version
for almost all~$t$, and then pointwise by the
continuity in~$t$.
It is readily seen from conditions (i) and (ii)
that the aforementioned properties, including (\ref{e2.3}),
also remain valid  for the functions
$$
w_\varepsilon(x)=\varepsilon ^{-d}g(x/\varepsilon),
$$
where $g$ is the standard Gaussian density and $\varepsilon\in (0,1)$.
Below we take for $\varepsilon$ only numbers of the
form~$1/n$, $n\in\mathbb{N}$.
Let us  set
$$
\varrho_\varepsilon:=\varrho * w_\varepsilon,
\quad
f_\varepsilon(t,x):=
\varrho_\varepsilon(t,x)+\varepsilon \max(1,|x|)^{-d-1},
$$
where we take for $\varrho_\varepsilon$ the version indicated
in~(\ref{e2.4}).
Since the function
$\varrho\Lambda$ is integrable,
one can find $\tau$ as close to $1$ as we like such that
\begin{equation}\label{e2.5}
\int_{\mathbb{R}^d} \varrho(\tau,x) \Lambda(x)\, dx<\infty.
\end{equation}
A number $\tau$ for which (\ref{e2.5}) is fulfilled can be
chosen in such a way that for each
$\varepsilon=1/n$, our version of $\varrho_\varepsilon(\tau,x)$
will coincide with the convolution
$\varrho(\tau,\,\cdot\,)*w_\varepsilon(x)$
for all~$x$.
Then the easily verified inequality
$\ln\max(|x+y|,1)\le \ln\max(|x|,1)+|y|$ yields
\begin{multline}\label{e2.6}
\int_{\mathbb{R}^d} f_\varepsilon(\tau,x)\Lambda(x)\, dx
\le
\int_{\mathbb{R}^d} \varrho_\varepsilon(\tau,x)\Lambda(x)\, dx
+
\varepsilon
\int_{\mathbb{R}^d} \max(|x|,1)^{-d-1}\Lambda(x)\, dx
\\
\le
\int_{\mathbb{R}^d} \varrho(\tau,x)\Lambda(x)\, dx
+
\int_{\mathbb{R}^d} |y|w_\varepsilon(y)\, dy
+
\varepsilon
\int_{\mathbb{R}^d} \max(|x|,1)^{-d-1}\Lambda(x)\, dx
\le M_1,
\end{multline}
where $M_1$ is a number independent of $\varepsilon$.
By (\ref{e2.3}) we have
\begin{equation}\label{e2.7}
\int_0^\tau\int_{\mathbb{R}^d}
\partial_t(\varrho * w_\varepsilon)\ln f_\varepsilon\, dx\, dt
=
\int_0^\tau\int_{\mathbb{R}^d}
\Bigl[
(a^{ij}\varrho)
*\partial_{x_i}\partial_{x_j}w_\varepsilon
-(b^i\varrho)*\partial_{x_i}w_\varepsilon\Bigr]\ln f_\varepsilon\,
dx\,dt,
\end{equation}
because $|f_\varepsilon|\le c_1+c_2\Lambda$
with some constants $c_1$ and $c_2$ and the functions
$((b^i\varrho)*\partial_{x_i}w_\varepsilon )\Lambda$
and
$((a^{ij}\varrho) *\partial_{x_i}\partial_{x_j}w_\varepsilon)
\Lambda$
are integrable on~$(0,1)\times\mathbb{R}^d$. Indeed,
since $\varrho *|\partial_{x_i}w_\varepsilon|>0$, one has
$$
((b^i\varrho)*\partial_{x_i}w_\varepsilon )\Lambda
=
((b^i\varrho)*\partial_{x_i}w_\varepsilon )
(|\partial_{x_i}w_\varepsilon|*\varrho)^{-1/2}
(|\partial_{x_i}w_\varepsilon|*\varrho)^{1/2}
\Lambda.
$$
Then
$((b^i\varrho)*\partial_{x_i}w_\varepsilon )
(|\partial_{x_i}w_\varepsilon|*\varrho)^{-1/2}
\in L^2((0,1)\times\mathbb{R}^d)$
by~(\ref{e2.2}) and the inclusion $|b|\in L^2(\mu)$.
In addition,
$(|\partial_{x_i}w_\varepsilon|*\varrho)^{1/2}\Lambda
\in L^2((0,1)\times\mathbb{R}^d)$
by the estimate
$$
|\ln\max(|x+y|,1)|^2\le 4+2|\ln\max(|x|,1)|^2+2|\ln\max(|y|,1)|^2
$$
and the same computations as in~(\ref{e2.6}).
Similarly we verify the integrability of the function
$[(a^{ij}\varrho) *\partial_{x_i}\partial_{x_j}w_\varepsilon]
\Lambda$ on~$(0,1)\times\mathbb{R}^d$.
We observe that one can integrate by parts on the right
in~(\ref{e2.6}).
Indeed,
$$
\int_0^\tau\int_{\mathbb{R}^d}
\frac{|\nabla f_\varepsilon|^2}{f_\varepsilon}\, dx\, dt
\le
2\int_0^\tau\int_{\mathbb{R}^d}
\frac{|\nabla \varrho_\varepsilon|^2}{\varrho_\varepsilon}\, dx\, dt
+2\varepsilon
(d+1)^2\int_{\{|x|\ge 1\}} |x|^{-d-3}\, dx,
$$
which is finite by (\ref{e2.2}), since
$\nabla \varrho_\varepsilon=(\nabla w_\varepsilon) *\varrho$
and $|\nabla w_\varepsilon |^2/w_\varepsilon \in L^1(\mathbb{R}^d)$.
In addition, one has
$[(b^i\varrho)*w_\varepsilon]\varrho_\varepsilon^{-1/2}\in
L^2((0,1)\times \mathbb{R}^d)$ again by
(\ref{e2.2}) and the inclusion $|b|\in L^2(\mu)$.
Similarly, we have
$[(a^{ij}\varrho)*\partial_{x_i}
w_\varepsilon]\varrho_\varepsilon^{-1/2}\in
L^2((0,1)\times \mathbb{R}^d)$.
Since $f_\varepsilon>\varrho_\varepsilon$, one has
$$
\partial_{x_i}f_\varepsilon\Bigl(
\partial_{x_j} [(a^{ij}\varrho)*w_\varepsilon]
-(b^i\varrho)* w_\varepsilon\Bigr) f_\varepsilon^{-1}\in
L^1(\mathbb{R}^d).
$$
Therefore,
\begin{equation}\label{e2.8}
\int_0^\tau\int_{\mathbb{R}^d}
\partial_t \varrho_\varepsilon  \ln f_\varepsilon\, dx\, dt
=-
\int_0^\tau\int_{\mathbb{R}^d}
\frac{\partial_{x_i}f_\varepsilon}{f_\varepsilon}
\Bigl(\partial_{x_j} [(a^{ij}\varrho)*w_\varepsilon]
-(b^i\varrho)* w_\varepsilon\Bigr)\, dx\, dt.
\end{equation}
The integrand on the left can be written as
$\partial_t(f_\varepsilon\ln f_\varepsilon)-\partial_t
\varrho_\varepsilon$.
Since the integrals of
$\varrho_\varepsilon(\tau,x)$ and $\varrho_\varepsilon(0,x)$ in
$x$ equal one,
we see that the left-hand side of (\ref{e2.8}) equals
$$
L_\varepsilon:=
\int_{\mathbb{R}^d}
[f_\varepsilon(\tau,x)\ln f_\varepsilon(\tau,x)-
f_\varepsilon(0,x)\ln f_\varepsilon(0,x)]\, dx .
$$
Since $f_\varepsilon(\tau,\,\cdot\,)
\ln f_\varepsilon(\tau,\,\cdot\,)\in L^1(\mathbb{R}^d)$
by (\ref{e2.6}) and the estimate
$|\ln f_\varepsilon|\le c_1+c_2\Lambda$, one has
$f_\varepsilon(0,\,\cdot\,)
\ln f_\varepsilon(0,\,\cdot\,)\in L^1(\mathbb{R}^d)$.
We need a lower bound on $L_\varepsilon$. To this end, we observe
that by the convexity of the function $s\mapsto s\ln s$ on
$(0,+\infty)$ and Jensen's inequality one has
\begin{multline*}
\int_{\mathbb{R}^d}
f_\varepsilon(0,x)\ln f_\varepsilon(0,x)\, dx
\\
\le
\int_{\mathbb{R}^d}
\varrho_\varepsilon(0,x)\ln (2\varrho_\varepsilon(0,x))\, dx
+
\int_{\mathbb{R}^d}
\varepsilon \max(|x|,1)^{-d-1}
\ln (2\varepsilon \max(|x|,1)^{-d-1})\, dx
\\
\le
\ln 2+
\int_{\mathbb{R}^d}
\varrho_\varepsilon(0,x)\ln \varrho_\varepsilon(0,x)\, dx
+
\varepsilon \ln 2\int_{\mathbb{R}^d} \max(|x|,1)^{-d-1}\, dx
\\
\le
\ln 2+
\int_{\mathbb{R}^d}
\varrho_0(x)\ln \varrho_0(x)\, dx
+
\varepsilon \ln 2\int_{\mathbb{R}^d} \max(|x|,1)^{-d-1}\, dx
=:M(\varepsilon).
\end{multline*}
On the other hand, (\ref{e2.6}) gives
$$
\int_{\mathbb{R}^d}
f_\varepsilon(\tau,x)\ln f_\varepsilon(\tau,x)\, dx
\ge
-(d+1)
\int_{\mathbb{R}^d}
f_\varepsilon(\tau,x)\Lambda(x)\, dx
-\varepsilon M_1(d+1)=:-K(\varepsilon).
$$
Note that for any bounded Borel
function $a$ on $(0,1)\times\mathbb{R}^d$
that
is Lipschitzian in the second argument
with Lipschitz norm~$\lambda$, for every $j$ we have
\begin{equation}\label{e2.9}
\partial_{x_j}[(a\varrho)*w_\varepsilon](t,x)=
a(t,x)\,\partial_{x_j}\,\varrho_\varepsilon(t,x)+
\int_{\mathbb{R}^d}
\partial_{x_j} w_\varepsilon(x-y)[a(t,y)-a(t,x)]\,\varrho(t,y)\, dy
\end{equation}
and
\begin{multline}\label{e2.10}
\Bigl|\int_{\mathbb{R}^d}
\partial_{x_j} w_\varepsilon(x-y)[a(t,y)-a(t,x)]\,\varrho(t,y)\,
dy\Bigr|
\le
\lambda\int_{\mathbb{R}^d}
|\partial_{x_j} w_\varepsilon(x-y)|\, |y-x|\,\varrho(t,y)\, dy
\\
\le
\lambda\int_{\mathbb{R}^d}
\varepsilon^{-d}\frac{|x-y|^2}{\varepsilon^2}\,
g\Bigl(\frac{x-y}{\varepsilon}\Bigr)\, \varrho(t,y)\, dy
=
\lambda(\varrho *q_\varepsilon)(t,x),
\end{multline}
where $q_\varepsilon(x):=w_\varepsilon(x)|x/\varepsilon|^2$,
$x\in \mathbb{R}^d$.
Let us note for the sequel that in the derivation
of (\ref{e2.8}) and (\ref{e2.9}) we have not used
the $\mu$-integrability of $\Lambda^2$ and the existence
of entropy of~$\mu_0$.
By using (\ref{e2.8}) and  (\ref{e2.9}) we obtain

\begin{multline*}
\int_0^1\int_{\mathbb{R}^d}
a^{ij}\frac{\partial_{x_i}f_\varepsilon}{f_\varepsilon}
\partial_{x_j}f_\varepsilon\, dx\, dt
=
\int_0^1\int_{\mathbb{R}^d}
\frac{\partial_{x_i}f_\varepsilon}{f_\varepsilon}
\Bigl[(b^i\varrho)* w_\varepsilon
+\varepsilon
a^{ij}\partial_{x_j}\frac{1}{\max(|x|,1)^{d+1}}
\Bigr]\, dx\, dt
\\
-\int_0^1\int_{\mathbb{R}^d}
\Bigl(\frac{\partial_{x_i}f_\varepsilon(t,x)}{f_\varepsilon(t,x)}
\int_{\mathbb{R}^d}
\partial_{x_j} w_\varepsilon(x-y)
[a^{ij}(t,y)-a^{ij}(t,x)]\,\varrho(t,y)\,
dy\Bigr)\, dx\, dt
-L_\varepsilon .
\end{multline*}
The right-hand side of this equality does not exceed
\begin{multline*}
\Bigl(\int_0^\tau\int_{\mathbb{R}^d}
\frac{|\nabla f_\varepsilon|^2}{f_\varepsilon}\, dx\, dt\Bigr)^{1/2}
\Bigl[\Bigl(\int_0^\tau\int_{\mathbb{R}^d}
\frac{\sum_{i=1}^d[(b^i\varrho)*w_\varepsilon]^2}{f_\varepsilon}\,
dx\, dt\Bigr)^{1/2}+\varepsilon MC_d\Bigr]
\\+
d^{3/2}\lambda
\Bigl(\int_0^\tau\int_{\mathbb{R}^d}
\frac{|\nabla f_\varepsilon|^2}{f_\varepsilon}\, dx\, dt\Bigr)^{1/2}
\Bigl(\int_0^\tau\int_{\mathbb{R}^d}
\frac{(\varrho*q_\varepsilon)^2}{f_\varepsilon}\, dx\, dt\Bigr)^{1/2}
+M(\varepsilon)+K(\varepsilon),
\end{multline*}
where
$M=\sup_{t,x}\|A(t,x)\|$
and $C_d$ is the integral of $(d+1)^2|x|^{-d-3}$ over the set
$\{|x|\ge 1\}$.
By (\ref{e2.2}) we have
$$
\int_0^\tau \int_{\mathbb{R}^d}
\frac{[(b^i\varrho)*w_\varepsilon]^2}{f_\varepsilon }\,
dx\, dt
\le  \int_0^\tau\int_{\mathbb{R}^d} |b^i|^2\, d\mu ,
\quad 1\le i\le d\ ,
$$
$$
\int_0^\tau\int_{\mathbb{R}^d}
\frac{(\varrho*q_\varepsilon)^2}{f_\varepsilon}\,
dx\, dt \le  \gamma:=\int_{\mathbb{R}^d} |x|^4 g(x)\, dx .
$$
Since $A\ge \alpha\cdot I$, we arrive at the estimate
\begin{multline*}
\alpha \int_0^\tau\int_{\mathbb{R}^d}
\frac{|\nabla f_\varepsilon|^2}{f_\varepsilon}
\, dx\, dt
\\
\le
\Bigl(\int_0^\tau \int_{\mathbb{R}^d}
\frac{|\nabla f_\varepsilon|^2}{f_\varepsilon}
\, dx\, dt\Bigr)^{1/2}
\Bigl(\|b\|_{2,\mu}+\varepsilon MC_d+\lambda d^{3/2}\sqrt{\gamma}\Bigr)
+
M(\varepsilon)+K(\varepsilon),
\end{multline*}
which by the inequality
$c\sqrt{x}\le \alpha x/2+c^2/(2\alpha)$ yields the estimate
\begin{equation}\label{e2.11}
\int_0^\tau\int_{\mathbb{R}^d}
\frac{|\nabla f_\varepsilon|^2}{f_\varepsilon}
\, dx\, dt
\le
\alpha^{-2}\Bigl(\|b\|_{2,\mu}+\varepsilon MC_d
+\lambda d^{3/2}\sqrt{\gamma}\Bigr)^2
+2\alpha^{-1}(M(\varepsilon)+K(\varepsilon)).
\end{equation}
The quantities $M(\varepsilon)$ and $K(\varepsilon)$ are uniformly
 bounded in~$\varepsilon$.
Letting $\varepsilon\to 0$ we obtain that
$\sqrt{\varrho(t,\,\cdot\,)}\in W^{2,1}(\mathbb{R}^d)$ for almost all $t\in (0,1)$.
Hence $\varrho(t,\,\cdot\,)\in W^{1,1}(\mathbb{R}^d)$
for almost all $t\in (0,\tau)$.
In addition,
the integral of $|\nabla \varrho|^2/\varrho$ does not exceed
the right-hand side of (\ref{e2.11}) with~$\varepsilon=0$.
Thus,
$\sqrt{\varrho}\in \mathbb{H}^{2,2}
([0,\tau]\times\mathbb{R}^d))$.
By the Sobolev embedding theorem we have
$\varrho\in L^{d/(d-2),1}([0,\tau]\times\mathbb{R}^d))$ if $d>2$, and
$\varrho\in L^{s,1}([0,\tau]\times\mathbb{R}^d))$
for all $s\in [1,\infty)$ if~$d=2$.

The last claim of the theorem is clear from our reasoning.
\end{proof}

The proof yields a useful estimate
\begin{multline}\label{e2.12}
\int_0^\tau \int_{\mathbb{R}^d}
\frac{|\nabla \varrho|^2}{\varrho}
\, dx\, dt
\le
\alpha^{-2}\Bigl(\|b\|_{2,\mu}+
\lambda d^{3/2}\sqrt{\gamma}\Bigr)^2
+2\ln 2 \alpha^{-1}
+
2\alpha^{-1}
\int_{\mathbb{R}^d}
\varrho_0(x)\ln \varrho_0(x)\, dx
\\
+
2\alpha^{-1}
(d+1)
\int_{\mathbb{R}^d} \varrho(\tau,x)\Lambda(x)\, dx.
\end{multline}

\begin{remark}
{\rm
It is clear from the proof that the entropy of
$\varrho_\varepsilon(0,x)$ has to be estimated only from above,
so in place of the integrability of
$\varrho(0,x)\ln \varrho(0,x)$
it suffices to require only the integrability of
$\varrho(0,x)\max(0,\ln \varrho(0,x))$
 (then Jensen's inequality must be applied to the function
$s\max(0,\ln s)$).
This leads to the effect that in estimate (\ref{e2.12}) in place
of $\varrho(0,x)\ln \varrho(0,x)$ we obtain
$\varrho(0,x)\max(0,\ln \varrho(0,x))$.
However, the obtained estimates and (\ref{e2.8})
show that  if we keep all other assumptions,
the entropy of $\varrho(0,x)$ is finite anyway.
But if no $\mu$-integrability of $\Lambda$ is required, then
the situation may change. For example, if $d=1$, $b=0$ and~$a=1/2$,
then for any initial distribution~$\mu_0$, the solution is given
by the convolution $\mu_0 *g_t$, where
$g_t(x)=(2\pi t)^{-1/2}\exp(-x^2/(2t))$.
If $\mu_0$ has a density $\varrho_0$ such that
$|\varrho_0'|^2/\varrho_0\in L^1(\mathbb{R}^1)$, but
the function $\varrho_0\ln\varrho_0$ is not integrable, then
the solution $\varrho(t,x)$ has no entropy for any~$t$,
although the quantities
 $\int |\partial_x \varrho(t,x)|^2\varrho(t,x)^{-1}\, dx$
 are uniformly bounded. The same example shows that
for validity of estimate (\ref{e2.1}) certain conditions on
the initial distribution are necessary.
It suffices to take for $\mu_0$ Dirac's measure at the origin.
Then the function $|\partial_x\varrho|^2/\varrho$ is not
integrable on $(0,1)\times\mathbb{R}^1$. It would be interesting
to find a sufficient condition on $A$ and $b$ ensuring
finite entropy of
 $\varrho(t,\,\cdot\,)$ for $t>0$ and any initial distribution.
}\end{remark}

In Example 3.1 below and in \cite{BDPR}
one can find conditions on the coefficients
$A$ and $b$ that ensure the inclusion $|b|\in L^2(\mu)$.

Estimate (\ref{e2.12}) can be improved under additional hypotheses
on $A$ and~$b$.

Set $b_0:=(b_0^j)$, $b_0^j=b^j-\partial_{x_i}a^{ij}$.

\begin{theorem}
Suppose
$\mu$ satisfies {\rm(\ref{e1.1})},~{\rm(\ref{e1.3})},
where $\nu=\varrho_0\, dx$, $\varrho_0$ has finite entropy and is
locally H\"older continuous.
Let $A$ and $b$ satisfy {\rm(C1)} and {\rm(C2)}  with some $p>d+2$.
 Suppose that
$|A^{-1/2}b_0|\in L^2(\mu)$,
$\ln(1+|x|)\in L^4(\mu)$
and that
\begin{equation}\label{e2.13}
\liminf\limits_{r\to\infty}
\int_0^1\int_{r\le |x|\le 2r}
\Bigl[r^{-4}\|A(t,x)\|^2
+ r^{-2}
\Theta_A(t,x)^2\Bigr]\, \mu_t(dx)\, dt=0.
\end{equation}
Then $\varrho(t,\,\cdot\,)\in W^{p,1}_{loc}$ and
for almost all $\tau\in [0,1]$ one has
$$
\int_0^\tau \int_{\mathbb{R}^d}
 \Bigl|\frac{\sqrt{A}\nabla \varrho}{\varrho}\Bigr|^2\, d\mu\le
\int_0^\tau\int_{\mathbb{R}^d}|A^{-1/2}b_0|^2\, d\mu
+2\int_{\mathbb{R}^d}[\varrho(0,x)\ln\varrho(0,x)
-\varrho(\tau,x)\ln\varrho(\tau,x)]\, dx
$$
and the right-hand side is finite.
Under the additional assumption that
$A\ge \alpha \cdot I$ for some $\alpha >0$, one has
$\sqrt{\varrho}\in \mathbb{H}^{2,2}([0,1]\times\mathbb{R}^d))$,
$\varrho\in L^{d/(d-2),1}([0,1]\times\mathbb{R}^d))$
if $d>2$ and
$\varrho\in L^{s,1}([0,1]\times\mathbb{R}^d))$
for all
$s\in [1,\infty)$ if~$d=2$.
\end{theorem}
\begin{proof}
By the local theory \cite{BKR01}, we know that
$\mu$ has a continuous positive density $\varrho$ such that
for every ball $B$ and every closed interval
$[t_1,t_2]$ in $(0,1)$ we have
$\|\varrho(t,\,\cdot\,)\|_{W^{p,1}(B)}\in L^p[t_1,t_2]$.
Let $B_j$ denote the closed
ball of radius $j$ centered at the origin.
We fix a function $\zeta\in C_0^\infty(\mathbb{R}^d)$ such that
$\zeta(x)=1$ if $|x|\le 1$, $\zeta(x)=0$ if
$|x|>2$. Set $\zeta_j(x):=\zeta(x/j)$.
For small $\varepsilon>0$ and large $k>0$,
let
$$
\varrho_{k,\varepsilon}:=\min(k,\varrho_\varepsilon),
\quad \varrho_k=\min(k,\varrho),
\quad
\Omega_{k,\varepsilon}:=\{\varrho_\varepsilon <k\},
\quad
\Omega_{k}:=\{\varrho <k\}.
$$
As in Theorem~\ref{t2.1}, for almost all $\tau$
one has~(\ref{e2.5}), which gives the integrability
of
$\varrho(\tau,\,\cdot\,)\ln \varrho(\tau,\,\cdot\,)$
on~$\mathbb{R}^d$.
For any $\delta> 0$ and $\tau=1-\delta$
we have the equality
\begin{multline*}
\int_\delta^\tau\int_{\mathbb{R}^d}
(\partial_t\varrho_\varepsilon)
(\ln\varrho_{k,\varepsilon})\zeta_j^2\,
dx\, dt
\\=
-\int_\delta^\tau
\int_{\mathbb{R}^d}
(a^{ik}\partial_{x_i}\varrho)*w_\varepsilon
\frac{\partial_{x_k}\varrho_\varepsilon}{\varrho_\varepsilon}
\zeta_j^2 I_{\Omega_{k,\varepsilon}}\,
dx\, dt
-2\int_\delta^\tau\int_{\mathbb{R}^d}
(a^{ik}\partial_{x_i}\varrho)*w_\varepsilon \partial_{x_k}\zeta_j
(\ln \varrho_{k,\varepsilon})\zeta_j\, dx\, dt
\\
+\int_\delta^\tau
\int_{\mathbb{R}^d} I_{\Omega_{k,\varepsilon}}
\Bigl((b_0\varrho)*w_\varepsilon,\frac{\nabla\varrho_\varepsilon}
{\varrho_\varepsilon}\Bigr)\, \zeta_j^2\, dx\, dt
+2\int_\delta^\tau\int_{\mathbb{R}^d}
\bigl((b_0\varrho)*w_\varepsilon,
\nabla\zeta_j)\bigr)\zeta_j\ln\varrho_{k,\varepsilon}\,
dx\, dt .
\end{multline*}
Since $\varrho_\varepsilon\partial_t\ln\varrho_{k,\varepsilon}
=\partial_t\varrho_{k,\varepsilon}$,
the left-hand side equals
\begin{multline*}
E(j,k,\varepsilon,\delta):=
\int_{\mathbb{R}^d} \zeta_j^2(x)\varrho_\varepsilon(\tau,x)
\ln\varrho_{k,\varepsilon}(\tau,x)\, dx
-\int_{\mathbb{R}^d} \zeta_j^2(x)\varrho_\varepsilon(\delta,x)
\ln\varrho_{k,\varepsilon}(\delta,x)\, dx
\\
+
\int_{\mathbb{R}^d} \zeta_j^2(x)
\varrho_{k,\varepsilon}(\delta,x)\, dx
-\int_{\mathbb{R}^d} \zeta_j^2(x)
\varrho_{k,\varepsilon}(\tau,x)\, dx.
\end{multline*}
Keeping  $\delta>0$ fixed, letting $\varepsilon\to 0$ and using
the integrability of the function
$t\mapsto \|\varrho(t,\,\cdot\,)\|_{W^{p,1}(B_{2j})}$
on $[\delta,\tau]$ as well as the
continuity and strict positivity of $\varrho$
on $[\delta,\tau]\times B_{2j}$, we obtain
\begin{multline*}
S_{j,k,\delta}:=
\int_\delta^\tau
\int_{\mathbb{R}^d}
\Bigl(A\nabla\varrho,
\frac{\nabla\varrho}{\varrho}\Bigr)
\zeta_j^2 I_{\Omega_{k}}\,
dx\, dt
\\
=
-2\int_\delta^\tau\int_{\mathbb{R}^d}
(A\nabla\varrho, \nabla \zeta_j)
\zeta_j\ln \varrho_{k}\, dx\, dt
+\int_\delta^\tau\int_{\mathbb{R}^d} I_{\Omega_{k}}
\Bigl(b_0,\frac{\nabla\varrho}
{\varrho}\Bigr)\, \zeta_j^2\varrho\, dx\, dt
\\
+2\int_\delta^\tau \int_{\mathbb{R}^d}
(b_0,\nabla\zeta_j)\zeta_j (\ln\varrho_{k})\,
\varrho\, dx\, dt
-E(j,k,\delta),
\end{multline*}
where
\begin{multline*}
E(j,k,\delta):=
\int_{\mathbb{R}^d} \zeta_j^2(x)\varrho(\tau,x)
\ln\varrho_{k}(\tau,x)\, dx
-\int_{\mathbb{R}^d} \zeta_j^2(x)\varrho(\delta,x)
\ln\varrho_{k}(\delta,x)\, dx
\\
+
\int_{\mathbb{R}^d} \zeta_j^2(x)
\varrho_{k}(\delta,x)\, dx
-\int_{\mathbb{R}^d} \zeta_j^2(x)
\varrho_{k}(\tau,x)\, dx.
\end{multline*}
Integrating by parts in the
integral of
$(A\nabla\varrho, \nabla \zeta_j)
\ln \varrho_{k}\zeta_j=
(\nabla\varrho, \zeta_jA\nabla \zeta_j)
\ln \varrho_{k}$
and writing
$( b_0,\nabla\varrho)=(A^{-1/2}b_0,
A^{1/2}\nabla\varrho)$,
we find
\begin{multline}\label{e2.14}
S_{j,k,\delta}
=2\int_\delta^\tau\int_{\Omega_{k}}
\Bigl(\frac{\nabla\varrho}{\varrho},A\nabla\zeta_j\Bigr)
\zeta_j\varrho\, dx \, dt
+2\int_\delta^\tau\int_{\mathbb{R}^d}
{\rm div}(\zeta_jA\nabla\zeta_j)
(\ln\varrho_{k})\, \varrho\, dx\, dt
\\
+\int_\delta^\tau\int_{\Omega_{k}}
\Bigl(b_0,\frac{\nabla\varrho}{\varrho}\Bigr)\,
\zeta_j^2\varrho\,
dx\, dt
+2\int_\delta^\tau\int_{\mathbb{R}^d}
(b_0,\nabla\zeta_j)\zeta_j (\ln\varrho_{k})\, \varrho\, dx\, dt
-E(j,k,\delta)\\
\le \sqrt{S_{j,k,\delta}}
\Bigl(2\| I_{\Omega_{k}}
\sqrt{A}\nabla \zeta_j\|_{L^2(\mu)}
+\| A^{-1/2}b_0\bigr\|_{L^2(\mu)}\Bigr)+R_{j,k,\delta}-E(j,k,\delta),
\end{multline}
where
$$
R_{j,k,\delta}:=2\int_\delta^\tau\int_{\mathbb{R}^d}
{\rm div}(\zeta_jA\nabla\zeta_j)
\ln\varrho_{k}\, \varrho\, dx\, dt
+2\int_\delta^\tau\int_{\mathbb{R}^d}
(b_0,\nabla\zeta_j)\zeta_j(\ln\varrho_{k})\, \varrho\, dx\, dt .
$$
Since $\varrho_0$ is H\"older continuous on~$B_{2j}$,
one has $\lim\limits_{\delta\to 0}
\varrho(\delta,x)=\varrho(0,x)$ uniformly on~$B_{2j}$
(see, e.g.,~\cite[Ch.~III, Theorem~7.1 and Theorem~10.1]{LSU}).
Hence
$$
\lim\limits_{\delta\to 0}
E(j,k,\delta)=E(j,k,0).
$$
Therefore, (\ref{e2.14}) holds for $\delta=0$.
Keeping~$k$ fixed, we observe that, given
$\varepsilon>0$,
for all sufficiently
large numbers $j$
of the form $j=r_l$ with $r_l\to\infty$
chosen according to (\ref{e2.13}),
the quantity
$R_{j,k,\delta}$ can
be made smaller than $\varepsilon$ in  absolute value.
Indeed, it follows by the hypotheses and the
estimates
$$
\sup_x |\nabla \zeta_j(x)|\le j^{-1}\sup_x |\nabla \zeta(x)|,
\quad
\sup_x |\partial_{x_i}\partial_{x_m}
\zeta_j(x)|\le j^{-2}
\sup_x |\partial_{x_i}\partial_{x_m} \zeta(x)|
$$
that for all $j=r_l$
the first term in the expression for $R_{j,k,\delta}$
can be estimated by
\begin{multline*}
M \|\ln \varrho_k\|_{L^2(\mu)} r_l^{-2}
\Bigl(\int_0^1\int_{\{r_l\le|x|\le 2r_l\}} \|A(t,x)\|^2\,
\mu_t(dx)\, dt\Bigr)^{1/2}
\\
+M \|\ln \varrho_k\|_{L^2(\mu)}
r_l^{-1}\Bigl(\int_0^1\int_{\{r_l\le|x|\le 2r_l\}} \Theta_A(t,x)^2\,
\mu_t(dx)\, dt\Bigr)^{1/2},
\end{multline*}
where $M$ is a constant that depends on the maxima of the
first and second derivatives of~$\zeta$.
The fact that $\ln\varrho_k\in L^2(\mu)$ follows
by the $\mu$-integrability of $|\ln (|x|+1)|^2$,
because on the set $\{x\colon\, \varrho(t,x)\le 1\}$ we have
$|\ln\varrho(t,x)|^2\sqrt{\varrho(t,x)}\le C$, hence
$$
|\ln\varrho(t,x)|^2\varrho(t,x)\le
(2d+2)^2(\ln(|x|+1))^2\varrho(t,x) +
C(|x|+1)^{-d-1}.
$$
Similarly, by the Cauchy inequality and the estimate
$$
|b(t,x)|\le \|A^{1/2}(t,x)\| |A^{-1/2}(t,x)b(t,x)|,
$$
the second term in the expression for
$R_{j,k,\delta}$ is majorized by
$$
Mr_l^{-1} \| A^{-1/2}b\|_{L^2(\mu)}
\|\ln\varrho_k\|_{L^4(\mu)}^2
\biggl(\int_0^1\int_{r_l\le |x|\le 2r_l} \|A(t,x)\|^2\,
\mu_t(dx)\, dt\biggr)^{1/2}.
$$
The quantities $E(j,k,0)$ are bounded from below by a constant
independent of $j$ and~$k$, because we consider only those
$\tau$ for which
$\varrho(\tau,\,\cdot\,)\ln \varrho(\tau,\,\cdot\,)
\in L^1(\mathbb{R}^d)$,
and we have
$\varrho(0,\,\cdot\,)\ln \varrho(0,\,\cdot\,)
\in L^1(\mathbb{R}^d)$ by assumption.
This yields that the integrals of
$|\sqrt{A}\nabla \varrho/\varrho|^2$ over the
sets $\Omega_{k}$
against~$\mu$ are uniformly bounded.
Letting $k\to\infty$ and then $j\to\infty$
we see that the function
$|\sqrt{A}\nabla \varrho|^2\varrho^{-1}$ is integrable on
 $[0,\tau]\times\mathbb{R}^d$. In addition, by (\ref{e2.14})
its integral $S$ satisfies the inequality
$S\le \sqrt{S}\| A^{-1/2}b_0\bigr\|_{L^2(\mu)}+ E$,
where $E$ is the difference of entropies of $\varrho(0,\,\cdot\,)$
and~$\varrho(\tau,\,\cdot\,)$.
This yields the desired bound.
\end{proof}

\begin{remark}
{\rm
 If $A$ is uniformly bounded, then the assumption
$\Lambda\in L^4(\mu)$ in the second theorem
can be relaxed to $\Lambda\in L^2(\mu)$.
}\end{remark}

\section{Higher integrability and boundedness of densities}

The results of the previous section show that the solutions
are globally integrable in some power greater than~$1$.
Here we derive yet stronger integrability properties and
the global boundedness under
additional assumptions on the coefficients.
In what follows we assume
that the measure $\mu$ is given by a density $\varrho$
such that for every $t\in [0,1)$, the function
$x\mapsto\varrho(t,x)$ is a probability density
with respect to Lebesgue measure.

Set
$$\|u\|_{p,q,\tau}=\|u I_{[0,\tau]}\|_{p,q},$$
where $t\to I_{[0,\tau]}(t)$ is the indicator function
of the interval $[0,\tau]$.

\begin{lemma}\label{lem3.1}
Let $d>2$. For every function $u\in
\mathbb{H}^{2,2}([0,\tau]\times\mathbb{R}^d)
\bigcap L^{2,\infty}([0,\tau]\times\mathbb{R}^d)$,
one has the inequality
$$\|u\|_{p,q,\tau}\le c(d, p)
\Bigl(\|\nabla u\|_{L^{2}([0,\tau]\times\mathbb{R}^d)}+
\|u\|_{2,\infty,\tau}\Bigr),
$$
where $2\le q, \, 2< p\le 2d/(d-2)$ and
$1/q+d/(2p)=d/4$.
\end{lemma}
\begin{proof}
Let $\delta=d/2 - d/p=2/q$.
Then
$$
\frac{(d-2)\delta}{2d}+\frac{1-\delta}{2}=\frac{1}{p}.
$$
Let $r=2d/(p\delta(d-2))$. Then
$r\ge 1$
($r>1$ if $p<2d/(d-2)$)  and
$r'=r/(r-1)=2/(p(1-\delta))$.
Writing $|u|^p=|u|^{p\delta}|u|^{p(1-\delta)}$ and
applying H\"older's inequality with $r$ and~$r'$
we obtain
$$
\Bigl(\int_{\mathbb{R}^d}
|u|^p\,dx\Bigr)^{1/p}\le
\Bigl(\int_{\mathbb{R}^d}|u|^{2d/(d-2)}\,dx\Bigr)
^{(d-2)\delta/d}
\Bigl(\int_{\mathbb{R}^d}|u|^2\,dx\Bigr)^{(1-\delta)/2}.
$$
The Sobolev inequality yields
$$
\|u(t,\,\cdot\,)\|_p\le c(d, p)
\|\nabla u(t,\,\cdot\,)\|_2^{\delta}
\, \|u(t,\,\cdot\,)\|_2^{1-\delta}.
$$
Therefore, for almost all $t\in [0,\tau]$ we have
$$
\|u(t,\cdot)\|_p\le c(d, p)
\|\nabla u(t,\,\cdot\,)\|_2^{\delta}
\, \|u\|_{2,\infty,\tau}^{1-\delta}.
$$
Since $\delta=2/q$, one has
$$
\Bigl(\int\limits_0^{\tau}
\|u(t,\,\cdot\,)\|_p^q\,dt\Bigr)^{1/q}
\le c(d, p)
\Bigl(\int\limits_0^{\tau}
\|\nabla u(t,\,\cdot\,)\|_2^2\,dt\Bigr)^{1/q}
\Bigl(\|u\|_{2,\infty,\tau}\Bigr)^{1-2/q}.
$$
By the Young inequality we obtain the required estimate.
\end{proof}

\begin{lemma}\label{lem3.2}
Suppose that hypothesis {\rm(i)}
 of Theorem {\rm\ref{t2.1}} is fulfilled
 and we have additionally
\begin{equation}\label{e3.1}
\sup_{t\in[0,1]}
\|b(t,\,\cdot\,)\|_{L^{s}(\mu_t)}<\infty,
\quad
\varrho\in L^{ks/(s-2)+1,k+(s-2)/s}([0,T]\times\mathbb{R}^d)
\end{equation}
with some $T\in [0,1]$, $s>2$, $k\ge 2/s$.
Furthermore, let
$\mu_0=\varrho(0,\,\cdot\,)\, dx$, where
$\varrho(0,\,\cdot\,)\in L^{k+1}(\mathbb{R}^d)$.
Then for almost all $\tau\in [0,T]$ we have
$\varrho(\tau,\,\cdot\,)\in W^{1,1}_{loc}$ and
\begin{multline}\label{e3.2}
\frac{2}{\alpha k(k+1)}\int_{\mathbb{R}^d}\varrho(\tau,x)^{k+1}\,dx+
\int_0^{\tau}\int_{\mathbb{R}^d}
|\nabla\varrho(t,x)|^2\varrho(t,x)^{k-1}\,dx\,dt
\\
\le
C(\alpha,\lambda,d,s)
\int_0^{\tau}\Bigl(\int_{\mathbb{R}^d}
\varrho(t,x)^{ks/(s-2)+1}\,dx\Bigr)^{(s-2)/s}\,dt
+\frac{2}{\alpha k(k+1)}
\int_{\mathbb{R}^d}\varrho(0,x)^{k+1}\,dx,
\end{multline}
where
$\alpha$ is the constant from the condition
$A(t,x)\ge \alpha \cdot I$ and
$C(\alpha,\lambda, d, s)$ is some number that depends
only on $\alpha, \lambda, d, s$.

If in place of {\rm(\ref{e3.1})} we have the condition
\begin{equation}\label{e3.3}
|b|\in L^{s}(\mu),
\quad
\varrho\in L^{ks/(s-2)+1}([0,T]\times\mathbb{R}^d),
\end{equation}
where $s>2$, $k>0$,
then for almost all $\tau\in [0,T]$, one has the inequality

\begin{multline}\label{e3.4}
\frac{2}{\alpha k(k+1)}
\int_{\mathbb{R}^d}
\varrho(\tau,x)^{k+1}\,dx+
\int_0^{\tau}\int_{\mathbb{R}^d}
|\nabla\varrho(t,x)|^2\varrho(t,x)^{k-1}\,dx\,dt
\\
\le
C(\alpha,\lambda,d,s)\Bigl(
\int_0^{\tau}\int_{\mathbb{R}^d}
\varrho(t,x)^{ks/(s-2)+1}
\,dx\,dt\Bigr)^{(s-2)/s}
+\frac{2}{\alpha k(k+1)}
\int_{\mathbb{R}^d}\varrho(0,x)^{k+1}\,dx.
\end{multline}
\end{lemma}
\begin{proof}
Let $\varrho_\varepsilon$ be the same as in Theorem~\ref{t2.1},
in particular, $\varepsilon\in\{1/n\}$.
By (\ref{e2.3}) we have
$$
\int_0^{\tau}\int_{\mathbb{R}^d}
\partial_{t}(\varrho*w_{\varepsilon})\varphi\,dx\,dt
=\int_0^{\tau}\int_{\mathbb{R}^d}
\Bigl[(a^{ij}\varrho)*\partial_{x_i}\partial_{x_j}
w_{\varepsilon}
-(b^i\varrho)*\partial_{x_i}w_{\varepsilon}\Bigr]
\varphi\,dx\,dt
$$
for every bounded measurable function $\varphi$
on~$[0,1]\times\mathbb{R}^d$,
because the indicated convolutions are integrable.
Let us take
$$
\varphi(t,x):=\varrho_{\varepsilon}(t,x)^k.
$$
We observe that
$$
\Bigl[(a^{ij}\varrho)*\partial_{x_j}w_{\varepsilon}
-(b^i\varrho)*w_{\varepsilon}\Bigr]
\varrho_{\varepsilon}^{k-1}\partial_{x_i}\varrho_{\varepsilon}
\in L^1([0,1]\times\mathbb{R}^d).
$$
Indeed, the function
$\varrho_{\varepsilon}^{k}$ is bounded and the functions
$$
\Bigl[(a^{ij}\varrho)*\partial_{x_j}w_{\varepsilon}
-(b^i\varrho)*w_{\varepsilon}\Bigr]
\varrho_{\varepsilon}^{-1/2},
\quad
|\nabla\varrho_\varepsilon| \varrho_{\varepsilon}^{-1/2}
$$
belong to $L^2([0,1]\times \mathbb{R}^d)$,
as already noted in the proof of Theorem~\ref{t2.1}.
Therefore, we can integrate by parts, which gives
\begin{equation}\label{e3.5}
\int_0^{\tau}\int_{\mathbb{R}^d}
\partial_{t}(\varrho*w_{\varepsilon})
\varrho^k_{\varepsilon}\,dx\,dt=
-\int_0^{\tau}\int_{\mathbb{R}^d}
\Bigl[(a^{ij}\varrho)*\partial_{x_j}w_{\varepsilon}-
(b^i\varrho)*w_{\varepsilon}\Bigr]k\varrho_{\varepsilon}^{k-1}
\partial_{x_i}\varrho_{\varepsilon}\,dx\,dt.
\end{equation}
Let $L_{\varepsilon}$ denote
 the left-hand side of this equality.
Then
$$
L_{\varepsilon}
=
\frac{1}{k+1}\int_{\mathbb{R}^d}
\varrho_{\varepsilon}(\tau,x)^{k+1}\,dx
-\frac{1}{k+1}\int_{\mathbb{R}^d}
\varrho_{\varepsilon}(0,x)^{k+1}\,dx .
$$
Note that the integrability of
$\varrho_{\varepsilon}(\tau,x)^{k+1}$ in $x$ follows by the
boundedness of this function and its integrability
for~$k=0$.
By using H\"older's inequality, we estimate
$L_{\varepsilon}$ from below as follows:
$$
L_{\varepsilon}\ge
\frac{1}{k+1}\int_{\mathbb{R}^d}
\varrho_{\varepsilon}(\tau,x)^{k+1}\,dx
-\frac{1}{k+1}\int_{\mathbb{R}^d}\varrho(0,x)^{k+1}\,dx.
$$
Let us consider the right-hand side  $R_{\varepsilon}$
of equality~(\ref{e3.5}).
By using equality~(\ref{e2.9}), we obtain
\begin{multline*}
R_{\varepsilon}=-k\int_0^{\tau}\int_{\mathbb{R}^d}
a^{ij}\partial_{x_j}
\varrho_{\varepsilon}\partial_{x_i}\varrho_{\varepsilon}
\varrho_{\varepsilon}^{k-1}\,dx\,dt
\\
-k\int_0^{\tau}\int_{\mathbb{R}^d}
\Bigl(\partial_{x_i}\varrho_{\varepsilon}(t,x)
\varrho_{\varepsilon}(t,x)^{k-1}
\int_{\mathbb{R}^d}
\partial_{x_j}w_{\varepsilon}(x-y)
[a^{ij}(t,y)-a^{ij}(t,x)]\varrho(t,y)\,dy\Bigr)\,dx\,dt
\\
+k\int_0^{\tau}\int_{\mathbb{R}^d}\varrho_{\varepsilon}
\partial_{x_i}\varrho_{\varepsilon}
\varrho^{k-1}_{\varepsilon}[(b^i\varrho)*w_{\varepsilon}]\,dx\,dt.
\end{multline*}
Hence, by (\ref{e2.10})
(we recall that in the derivation of (\ref{e2.10}) we have not
used the $\mu$-integrability of $|\ln (1+|x|)|^2$ and the existence
of entropy of~$\mu_0$)
we have
\begin{multline*}
R_{\varepsilon}\le -k\alpha\int_0^{\tau}\int_{\mathbb{R}^d}
|\nabla\varrho_{\varepsilon}|^2
\varrho^{k-1}_{\varepsilon}\,dx\,dt
\\
+kd^{3/2}\lambda\Bigl(\int_0^{\tau}\int_{\mathbb{R}^d}
|\nabla\varrho_{\varepsilon}|^2
\varrho^{k-1}_{\varepsilon}\,dx\,dt\Bigr)^{1/2}
\Bigl(
\int_0^{\tau}\int_{\mathbb{R}^d}\varrho_{\varepsilon}^{k-1}
(\varrho*q_{\varepsilon})^2\,dx\,dt
\Bigr)^{1/2}
\\
+k\Bigl(
\int_0^{\tau}\int_{\mathbb{R}^d}|\nabla\varrho_{\varepsilon}|^2
\varrho^{k-1}_{\varepsilon}\,dx\,dt\Bigr)^{1/2}
\Bigl(\int_0^{\tau}\int_{\mathbb{R}^d}\varrho^{k-1}_{\varepsilon}
\sum_{i=1}^d[(b^i\varrho)*w_{\varepsilon}]^2
\,dx\,dt\Bigr)^{1/2},
\end{multline*}
where $q_{\varepsilon}$ is the same as in~(\ref{e2.10}).
By the inequality $ab\le \alpha a^2/4 + b^2/\alpha$
we obtain the estimate
\begin{multline*}
R_{\varepsilon}\le -\frac{1}{2}k\alpha
\int_0^{\tau}\int_{\mathbb{R}^d}
|\nabla\varrho_{\varepsilon}|^2
\varrho^{k-1}_{\varepsilon}\,dx\,dt
+\frac{kd^3\lambda^2}{\alpha}
\int_0^{\tau}\int_{\mathbb{R}^d}\varrho_{\varepsilon}^{k-1}
(\varrho*q_{\varepsilon})^2\,dx\,dt
\\
+\frac{k}{\alpha}
\int_0^{\tau}\int_{\mathbb{R}^d}\varrho^{k-1}_{\varepsilon}
\sum_{i=1}^d[(b^i\varrho)*w_{\varepsilon}]^2
\,dx\,dt.
\end{multline*}
Combining our bounds on $L_\varepsilon$ and $R_\varepsilon$
we arrive at the inequality
\begin{multline*}
\frac{2}{\alpha k(k+1)}\int_{\mathbb{R}^d}
\varrho_{\varepsilon}(\tau,x)^{k+1}\,dx
+
\int_0^{\tau}\int_{\mathbb{R}^d}
|\nabla\varrho_{\varepsilon}|^2
\varrho^{k-1}_{\varepsilon}\,dx\,dt
\\
\le C(\alpha,\lambda,d)\Bigl(
\int_0^{\tau}\int_{\mathbb{R}^d}
\varrho_{\varepsilon}^{k-1}(\varrho*q_{\varepsilon})^2\,dx\,dt
+
\int_0^{\tau}\int_{\mathbb{R}^d}\varrho^{k-1}_{\varepsilon}
\sum_{i=1}^d[(b^i\varrho)*w_{\varepsilon}]^2\,dx\,dt\Bigr)
\\
+\frac{2}{\alpha k(k+1)}
\int_{\mathbb{R}^d}\varrho(0,x)^{k+1}\,dx.
\end{multline*}
Note that
$|(b^i\varrho)* w_\varepsilon|^2\le
\varrho_\varepsilon (|b^i|^2\varrho) *w_\varepsilon$.
By H\"older's inequality we have
\begin{multline*}
\int_{\mathbb{R}^d}\varrho^{k-1}_{\varepsilon}
|(b^i\varrho)*w_{\varepsilon}|^2(t,x)\,dx
\\
\le \int_{\mathbb{R}^d}\varrho^{k}_{\varepsilon}
(|b^i|^2\varrho)*w_{\varepsilon}(t,x)\,dx=
\int_{\mathbb{R}^d}
\Bigl(\int_{\mathbb{R}^d}
\varrho_{\varepsilon}(t,x)^{k}
w_{\varepsilon}(x-y)\,dx\Bigr)
|b^i(t,y)|^2\varrho(t,y)\,dy
\\
\le
\|(b^i(t,\,\cdot\, ))^2\|
_{L^{s/2}(\mu_t)}\Bigl(\int_{\mathbb{R}^d}
\Bigl(
\int_{\mathbb{R}^d}
\varrho_{\varepsilon}(t,x)^kw_{\varepsilon}(x-y)\,dx
\Bigr)^{s/(s-2)}\varrho(t,y)\,dy\Bigr)^{(s-2)/s}
\\
\le
\|b^i(t,\,\cdot\,)\|^2_{L^s(\mu_t)}
\Bigl(
\int_{\mathbb{R}^d}
\varrho_{\varepsilon}(t,x)^{k(s-2)/s}
w_{\varepsilon}(x-y)\,dx
\Bigr)\varrho(t,y)\,dy\Bigr)^{(s-2)/s}
\\
=
\|b^i(t,\,\cdot\,)\|^2_{L^s(\mu_t)}
\Bigl(\int_{\mathbb{R}^d}
\varrho_{\varepsilon}(t,x)^{ks/(s-2)+1}
\,dx\Bigr)^{(s-2)/s}.
\end{multline*}
Similarly, taking into account the estimate
$$
|\varrho *q_\varepsilon(t,x)|^2\le
\varrho_\varepsilon(t,x)
\int_{\mathbb{R}^d}
\varrho(t,y)\frac{|x-y|^4}{\varepsilon^4}w_\varepsilon(x-y)\, dy
$$
we obtain
\begin{multline*}
\int_{\mathbb{R}^d}\varrho_{\varepsilon}^{k-1}
(\varrho*q_{\varepsilon})^2(t,x)\,dx
\le
\int_{\mathbb{R}^d}
\int_{\mathbb{R}^d}
\varrho_{\varepsilon}(t,x)^{k}
\varrho(t,y)\frac{|x-y|^4}{\varepsilon^4}w_\varepsilon(x-y)\, dx\, dy
\\
\le
\gamma(s)\int_{\mathbb{R}^d}
\varrho(t,y)\Bigl(\int_{\mathbb{R}^d}
\varrho_\varepsilon(t,x)^{ks/(s-2)}
w_\varepsilon(x-y)\,
dx\Bigr)^{s/(s-2)}\, dy
\\
\le
\gamma(s)\Bigl(\int_{\mathbb{R}^d}\int_{\mathbb{R}^d}
\varrho(t,y)
\varrho_\varepsilon(t,x)^{ks/(s-2)}
w_\varepsilon(x-y)\,
dx\, dy \Bigr)^{s/(s-2)}
\\
=
\gamma(s)
\Bigl(\int_{\mathbb{R}^d}
\varrho_{\varepsilon}(t,x)^{ks/(s-2)+1}
\,dx\Bigr)^{(s-2)/s},
\end{multline*}
where
$$
\gamma(s):=\Bigl(
\int_{\mathbb{R}^n}|x|^{2s}g(x)\,dx\Bigr)^{2/s}.
$$
Set
$$
B(s):=\max_{t\in[0,1]}
\bigl\| b(t,\,\cdot\, )\bigr\|_{L^s(\mu_t)}^2.
$$
Finally, we obtain
\begin{multline*}
\frac{2}{\alpha k(k+1)}\int_{\mathbb{R}^d}
\varrho_{\varepsilon}(\tau,x)^{k+1}\,dx
+
\int_0^{\tau}\int_{\mathbb{R}^d}
|\nabla\varrho_{\varepsilon}|^2
\varrho^{k-1}_{\varepsilon}\,dx\,dt
\\
\le
C(\alpha,\lambda,d)\Bigl(\gamma(s)+B(s)\Bigr)
\int_0^{\tau}\Bigl(\int_{\mathbb{R}^d}
\varrho_{\varepsilon}(t,x)^{ks/(s-2)+1}
\,dx\Bigr)^{(s-2)/s}\,dt
\\
+
\frac{2}{\alpha k(k+1)}\int_{\mathbb{R}^d}\varrho(0,x)^{k+1}\,dx.
\end{multline*}
Passing to the limit as $\varepsilon\to 0$, we obtain
the required estimate for almost every $\tau\in [0,T]$.
Indeed, for almost every~$\tau$ we have
$\lim\limits_{\varepsilon\to 0}
\varrho_\varepsilon(\tau,x)=\varrho(\tau,x)$ for
almost all~$x$,
and the right-hand side of the above inequality remains
bounded as $\varepsilon\to 0$.
This yields that $\varrho(\tau,\,\cdot\,)\in W^{1,1}_{loc}$
for almost all~$\tau$. Hence
$\lim\limits_{\varepsilon\to 0}
\nabla \varrho_\varepsilon(\tau,x)=\nabla \varrho(\tau,x)$ for
almost all~$x$, which gives the indicated estimate
by Fatou's theorem.

In the case of condition (\ref{e3.3})
the proof  repeats entirely the foregoing reasoning except for
the estimate of the integral of
$\varrho^{k-1}_{\varepsilon}|(b^i\varrho)*w_{\varepsilon}|^2$,
for which we have
\begin{multline*}
\int_{0}^{\tau}\int_{\mathbb{R}^d}\varrho^{k-1}_{\varepsilon}
|(b^i\varrho)*w_{\varepsilon}|^2\,dx\,dt
\\
\le \int_{0}^{\tau}\int_{\mathbb{R}^d}\varrho^{k}_{\varepsilon}
(|b^i|^2\varrho)*w_{\varepsilon}\,dx\,dt=
\int_{0}^{\tau}\int_{\mathbb{R}^d}
\Bigl(\int_{\mathbb{R}^d}\varrho_{\varepsilon}(t,x)^k
w_{\varepsilon}(x-y)\,dx\Bigr)
|b^i(t,y)|^2\varrho(t,y)\,dy\,dt
\\
\le
\|(b^i)^2\|_{L^{s/2}(\mu)}\Bigl(\int_{0}^{\tau}
\int_{\mathbb{R}^d}
\Bigl(\int_{\mathbb{R}^d}
\varrho_{\varepsilon}(t,x)^k w_{\varepsilon}(x-y)\,dx
\Bigr)^{s/(s-2)}\varrho(t,y)\,dy\,dt\Bigr)^{(s-2)/s}
\\
\le
\|b\|_{L^s(\mu)}^2
\Bigl(\int_{0}^{\tau}\int_{\mathbb{R}^d}
\varrho_{\varepsilon}(t,x)^{ks/(s-2)+1}\,dx\,dt
\Bigr)^{(s-2)/s}
\end{multline*}
by H\"older's inequality.
\end{proof}

\begin{remark}\label{rem3.1}
{\rm
It is seen from the proof that the assumption that the integrals
of $\varrho(t,x)$ with respect to $x$ equal $1$ can be replaced
by the assumption that these integrals are uniformly bounded.
}\end{remark}

\begin{theorem}\label{t3.1}
Suppose that under the hypotheses of Theorem {\rm\ref{t2.1}}
we have additionally
$$
\sup_{t\in[0,1]}
\|b(t,\,\cdot\,)\|_{L^d(\mu_t)}<\infty
$$
and $\mu_0=\varrho(0,\,\cdot\,)\, dx$, where
$\varrho(0,\,\cdot\,)\in L^p(\mathbb{R}^d)$ for all
$p\in [1,+\infty)$.
Then
$$
\varrho\in L^{p,q}([0,\tau]\times \mathbb{R}^d)
$$
for all $p,q\in [1,+\infty)$ and $\tau\in (0,1)$.
\end{theorem}
\begin{proof}
Let us consider the case $d>2$.
Let $\varrho\in
L^{ks/(s-2)+1, k+(s-2)/s}([0,\tau]\times\mathbb{R}^d)$,
where $s>2$ and~$k\ge 2/s$. Then
by Lemma \ref{lem3.2} we obtain the estimate
\begin{multline*}
\int^{\tau}_0\int_{\mathbb{R}^d}|\nabla\varrho(t,x)|^2\varrho(t,x)^{k-1}\,dx\,dt
\\
\le
C(\alpha,\lambda,d,s)
\int_0^{\tau}\Bigl(\int_{\mathbb{R}^d}
\varrho(t,x)^{ks/(s-2)+1}\,dx\Bigr)^{(s-2)/s}\,dt
+\frac{2}{\alpha k(k+1)}\int_{\mathbb{R}^d}\varrho(0,x)^{k+1}\,dx.
\end{multline*}
By the Sobolev   embedding theorem applied to the functions
$x\mapsto \varrho(t,x)^{(k+1)/2}$ with
$$
|\nabla \varrho^{(k+1)/2}|^2
=\frac{(k+1)^2}{4}|\varrho^{(k-1)/2}\nabla \varrho|^2
$$
we have
\begin{multline*}
\frac{4}{(k+1)^2}\int_0^{\tau}\Bigl(\int_{\mathbb{R}^d}
\varrho(t,x)^{d(k+1)/(d-2)}\,dx\Bigr)^{(d-2)/d}\,dt
\le
\int_0^{\tau}\int_{\mathbb{R}^d}
|\nabla\varrho(t,x)|^2
\varrho(t,x)^{k-1}\,dx\,dt
\\
\le \widetilde{C}(\alpha,\lambda,d,s)
\int_0^{\tau}\Bigl(\int_{\mathbb{R}^d}
\varrho(t,x)^{ks/(s-2)+1}
\,dx\Bigr)^{(s-2)/s}\,dt
+\frac{2}{\alpha k(k+1)}
\int_{\mathbb{R}^d}\varrho(0,x)^{k+1}\,dx.
\end{multline*}
Thus, one has
\begin{multline}\label{e3.6}
\Bigl(\|\varrho\|_{d(k+1)/(d-2),k+1,\tau}\Bigr)^{k+1}
\\
\le
C(\alpha,\lambda,d,s,k)
\Bigl(\|\varrho\|_{ks/(s-2)+1,k+(s-2)/s,\tau}
\Bigr)^{k+(s-2)/s}
+M(k,\alpha),
\end{multline}
where
$$
M(k,\alpha):=\frac{k+1}{2\alpha k}
\int_{\mathbb{R}^d}\varrho(0,x)^{k+1}\,dx .
$$
Now set
$$
p_n:=p_{n-1}+\frac{2}{d-2},
\quad
q_n:=q_{n-1}+\frac{2}{d},
\quad
q_1=1, \quad p_1=\frac{d}{d-2}.
$$
By Theorem  \ref{t2.1} we have
$\varrho\in L^{p_1,1}([0,\tau]\times\mathbb{R}^d)$
for all~$\tau<1$.
This enables us to start iterations based on~(\ref{e3.6}).
Namely, if in (\ref{e3.6}) we set
$s=d$ and $k=q_{n-1}-(d-2)/d=q_n-1$,
then we arrive at the estimate
$$
\Bigl(\|\varrho\|_{p_n,q_n,\tau}\Bigr)^{q_n}\le
C(\alpha,\lambda,d,d,q_n-1)
\Bigl(\|\varrho\|_{p_{n-1},q_{n-1},\tau}\Bigr)^{q_{n-1}}
+M(q_{n}-1,\alpha).
$$
Since $p_n\to\infty$ and $q_n\to\infty$
as $n\to\infty$, the theorem is proven in the case $d>2$.
The cases $d=1$ and $d=2$ are even simpler, because in the Sobolev inequality
in place of the exponent $d/(d-2)$ one can take any number~$r>1$. However,
we need not consider these cases separately  and deduce them from the result for $d=3$.
To this end, we pass from the function of two variables to the function
of three variables
$u=\varrho(t,x_1,x_2)g(x_3)$, where $g$ is the standard Gaussian density.
The measure $u\, dx\, dt$ satisfies our equation on
$[0,1)\times \mathbb{R}^3$
with the coefficients $a^{ij}$ and $b^i$ that
coincide with the initial ones if $i,j\le 2$, and
$a^{33}=1$, $a^{3j}=a^{j3}=0$, $b^3(t,x)=-x_3$.
\end{proof}

\begin{theorem}\label{t3.2}
Suppose that under the hypotheses of Theorem {\rm\ref{t2.1}}
for some $\beta>d+2$ we have
$|b|\in L^\beta(\mu)$
and $\varrho(0,\,\cdot\,)\in L^{\infty}(\mathbb{R}^d)$.
Supppose that either
$\sup_{t\in[0,1]}
\|b(t,\,\cdot\,)\|_{L^d(\mu_t)}<\infty
$
or
$\varrho\in L^{p}([0,\tau]\times\mathbb{R}^d)$,
for all $\tau<1$ with some $p>1$.
Then $\varrho\in L^{\infty}([0,\tau]\times \mathbb{R}^d)$
for every $\tau<1$.
\end{theorem}
\begin{proof}
Let $d>2$.
Let us fix $\tau<1$. Let
$\varrho\in
L^{k\beta/(\beta -2)+1}
([0,\tau]\times\mathbb{R}^d)$, where $k>0$.
Then by (\ref{e3.4}) we obtain
\begin{multline*}
\int_0^{\tau}
\int_{\mathbb{R}^d}|\nabla\varrho(t,x)|^2
\varrho(t,x)^{k-1}\,dx\,dt
\\
\le
C(\alpha,\lambda,d,\beta)\Bigl(
\int_0^{\tau}\int_{\mathbb{R}^d}
\varrho(t,x)^{k\beta/(\beta -2)+1}
\,dx\,dt\Bigr)^{(\beta -2)/\beta}
+\frac{2}{\alpha k(k+1)}
\int_{\mathbb{R}^d}\varrho(0,x)^{k+1}\,dx.
\end{multline*}
Therefore,
\begin{multline*}
\frac{4}{(k+1)^2}\Bigl(\|\nabla
(\varrho^{(k+1)/2})\|_{L^{2}([0,\tau]\times\mathbb{R}^d)}
\Bigr)^2
\\
\le
C(\alpha,\lambda,d,\beta)\Bigl(
\int_0^{\tau}\int_{\mathbb{R}^d}
\varrho(t,x)^{k\beta/(\beta -2)+1}
\,dx\,dt\Bigr)^{(\beta -2)/\beta}
+\frac{2}{\alpha k(k+1)}
\int_{\mathbb{R}^d}\varrho(0,x)^{k+1}\,dx.
\end{multline*}
By inequality \ref{e3.4}  we have for almost all $t<\tau$
\begin{multline*}
\int_{\mathbb{R}^d}\varrho(t,x)^{k+1}\,dx
\\
\le
k(k+1)C(\alpha,\lambda,d,\beta)\Bigl(
\int_0^{\tau}\int_{\mathbb{R}^d}
\varrho^{k\beta/(\beta -2)+1}\,dx\,dt
\Bigr)^{(\beta -2)/\beta}
+\int_{\mathbb{R}^d}\varrho(0,x)^{k+1}\,dx,
\end{multline*}
whence we obtain
\begin{multline*}
\Bigl(\|\varrho^{(k+1)/2}\|_{2,\infty,\tau}\Bigr)^2
\\
\le
k(k+1)C(\alpha,\lambda,d,\beta)\Bigl(
\int_0^{\tau}\int_{\mathbb{R}^d}
\varrho^{k\beta/(\beta -2)+1}\,dx\,dt
\Bigr)^{(\beta -2)/\beta}
+\int_{\mathbb{R}^d}\varrho(0,x)^{k+1}\,dx.
\end{multline*}
Lemma \ref{lem3.1} for
$u=\varrho^{(k+1)/2}$ and
$p=q=2(d+2)/d$
 yields the estimate
$$
\|\varrho^{(k+1)/2}\|_{p,q,\tau}\le
C(d,2(d+2)/d)\Bigl(
\|\nabla(
\varrho^{(k+1)/2})\|_{L^2([0,\tau]\times \mathbb{R}^d)}
+\|\varrho^{(k+1)/2}\|_{2,\infty,\tau}\Bigr).
$$
By using the inequality $(a+b)^2\le 2a^2+2b^2$ and
the bounds found above we obtain
\begin{multline}\label{e3.7}
\Bigl(
\|\varrho\|_{L^{(d+2)(k+1)/d}([0,\tau]\times\mathbb{R}^d)}
\Bigr)^{k+1}
\le 4(k+1)^2C(\alpha,\lambda,d,\beta)\Bigl(
\|\varrho\|_{L^{k\beta/(\beta-2)+1}
([0,\tau]\times\mathbb{R}^d)}
\Bigr)^{k+(\beta -2)/\beta}
\\
+\Bigl(2+\frac{k+1}{\alpha k}\Bigr)
\int_{\mathbb{R}^d}\varrho(0,x)^{k+1}\,dx.
\end{multline}
Since $\varrho(0,\,\cdot\,)\in
L^{\infty}([0,\tau]\times\mathbb{R}^d)$, there exists
a constant $C$ such that
$$
\Bigl(2+\frac{k+1}{\alpha k}\Bigr)
\int_{\mathbb{R}^d}\varrho(0,x)^{k+1}\,dx
\le C^{k+1}.
$$
Now we set
$$
p_n=\frac{(d+2)(\beta-2)}{d\beta}
\Bigl(p_{n-1}+\frac{2}{\beta-2}\Bigr),
\quad
p_1=\frac{d}{d-2},
$$
$$
A_n=\|\varrho\|_{L^{p_n}([0,\tau]\times\mathbb{R}^d)},
\quad
C_1=4\Bigl(\frac{d}{d+2}\Bigr)^2C(\alpha,\lambda, d,\beta).
$$
We have $\varrho\in
L^{p_1}([0,\tau]\times\mathbb{R}^d)$
by Theorem \ref{t3.1}
if $\sup_t
\|b(t,\,\cdot\,)\|_{L^d(\mu_t)}<\infty$;
if $\varrho\in L^p$ с $p>1$, then the same can be readily
deduced from~(\ref{e3.7}).
Note that
$$
1< \frac{(d+2)(\beta-2)}{d\beta}\le \frac{d}{d-2},
$$
since $\beta>d+2$.
Hence $p_n\ge \bigl((d+2)(\beta-2)/(d\beta)\bigr)^{n}$.
In order to prove the membership of $\varrho$ in the
space $L^\infty([0,\tau]\times\mathbb{R}^d)$,
it suffices to establish the uniform boundedness of~$\{A_n\}$.
Suppose that
 $A_n\to\infty$. Then there exists $N$ such
 that for every $n>N$ we have
$$
\frac{C}{A_{n-1}}<1,\quad
 (A_{n-1})^{-2/\beta}<1,
 \quad
  \Bigl(\frac{C}{A_{n-1}}\Bigr)^{p_nd/(d+2)}<1.
 $$
Let $k:=p_nd/(d+2)-1$. Then
$$
p_{n-1}=\frac{k\beta}{\beta -2}+1,
\quad
p_{n-1}\frac{\beta -2}{\beta} -p_n\frac{d}{d+2}=-\frac{2}{\beta}
$$
so (\ref{e3.7}) yields
$$
(A_n)^{p_nd/(d+2)}\le
p_n^2 C_1(A_{n-1})^{p_{n-1}(\beta -2)/\beta}
+C^{p_nd/(d+2)}.
$$
Therefore,
$$
\Bigl(\frac{A_n}{A_{n-1}}\Bigr)^{p_nd/(d+2)}\le
p_n^2 C_1(A_{n-1})^{-2/\beta}+
\Bigl(\frac{C}{A_{n-1}}\Bigr)^{p_nd/(d+2)}
\le p_n^2C_1+1.
$$
whence we obtain
$$
\ln A_n-\ln A_{n-1}\le
\frac{d+2}{p_nd}\ln\Bigl(p_n^2C_1+1\Bigr), \quad n>N.
$$
Since $p_n\ge \bigl((d+2)(\beta-2)/(d\beta)\bigr)^{n}$,
the above estimate yields
convergence of the series of $\ln A_n-\ln A_{n-1}$, which
contradicts the unboundedness of~$A_n$. The case $d\le 2$ is justified in the same
manner as in the previous theorem.
\end{proof}

\begin{remark}\label{rem3.2}
{\rm
(i)
According to Remark 2.1, no integrability
$\varrho(0,x)\ln\varrho(0,x)$ in Theorems 3.1 and 3.2 is required,
since the integrability of
$\varrho(0,x)\max(0,\ln\varrho(0,x))$ follows by the inclusion
$\varrho(0,\,\cdot\,)\in L^p(\mathbb{R}^d)$ with~$p>1$.

(ii)
It is seen from the proof and Remark \ref{rem3.1}
that the assumption in Theorem \ref{t3.1} and Theorem \ref{t3.2}
 that the integrals
of $\varrho(t,x)$ with respect to $x$ equal $1$ can be replaced
by the assumption that these integrals are uniformly bounded.

(iii) Let us note that if in Theorem \ref{t3.2}
it is given in advance that
$\varrho\in L^{p}([0,\tau]\times\mathbb{R}^d)$
for some $p>1$, then we need not require the integrability
of the function $|\ln(1+|x|)|^2\varrho(t,x)$, but the boundedness of
$\varrho(0,x)$ is important.
}\end{remark}

Now we employ the proven theorem for obtaining
upper bounds on~$\varrho$.
As in the elliptic case considered in the papers
 \cite{MPR}, \cite{BKR05a}, \cite{BKR05},
the idea is this: in order to obtain a pointwise
estimate $\varrho(t,x)\le \Phi(t,x)^{-1}$,
one has to consider
the measure $\nu:=\Phi\cdot\mu$ and establish the boundedness
of its density.
We shall consider functions $\Phi$ that do not
depend on~$t$.
If $\Phi$ has locally bounded first and second order
derivatives, then the measure $\nu$ satisfies the
equation
$$
L^{*}\nu=(a^{ij}\partial_{x_i}\partial_{x_j}\Phi)\varrho
+2\partial_{x_i}\Phi\partial_{x_j}(a^{ij}\varrho)
-b^i \partial_{x_i}\Phi \varrho
=-L\Phi\cdot\varrho
+ 2\partial_{x_j}(a^{ij}\partial_{x_i}\Phi\varrho )
$$
understood in the same sense as~(\ref{e1.1}).

\begin{theorem}\label{t3.3}
Suppose that all hypotheses of Theorem~{\rm\ref{t3.2}}
are fulfilled and we are given a function
 $\Phi\ge c>0$ on $\mathbb{R}^d$ with
locally bounded second order derivatives such that
$\varrho(0,x)\le C\Phi(x)^{-1}$, $\Phi\in L^1(\mu_0)$
and
$$
\Phi^{1+\varepsilon},\
 |L\Phi|^{\beta/2}\Phi^{1-\beta/2},\
  |A\nabla\Phi|^\beta \Phi^{1-\beta}\in L^1(\mu),
\quad
\sup_{t\in [0,1]}
\int_{\mathbb{R}^d} \Phi(x)\varrho(t,x)\, dx <\infty
$$
with some $\varepsilon>0$.
Then for every $\tau<1$ there is a number $C_\tau$
such that
$$
\varrho(t,x)\le C_\tau \Phi(x)^{-1}
\quad
\hbox{for almost all $(t,x)\in [0,\tau)$.}
$$
\end{theorem}
\begin{proof}
It is seen from the reasoning used in the proof of
Theorem \ref{t3.2} and Remark \ref{rem3.2}
that it suffices to establish
estimate (\ref{e3.4}) with $s=\beta$
for the measure $\nu=\Phi\cdot\mu$ whose density
belongs to $L^{1+\varepsilon}([0,1]\times \mathbb{R}^d)$
by the boundedness of~$\varrho$.
In that estimate a homogeneous equation was concerned,
and the measures $\mu_t$ were probabilities.
However, under present assumptions the same
estimate remains valid  in the presence of the indicated
right-hand side as well if in place of the condition
$\mu_t(\mathbb{R}^d)=1$ we assume only the uniform
boundedness of measures~$\mu_t$.
Indeed, on the right-hand side of
(\ref{e3.5}) with $\nu$ in place of~$\mu$,
i.e., with $u:=\Phi\varrho$ in place of~$\varrho$,
there appears additionally the integral of the
expression
$$
-(L\Phi\cdot\varrho)
*w_\varepsilon [u*w_\varepsilon]^k
+ [2\partial_{x_j}(a^{ij}\partial_{x_i}\Phi\varrho)]*w_\varepsilon
[u*w_\varepsilon]^k .
$$
Let us set $\xi:=L\Phi/\Phi$, $\eta:=|A\nabla\Phi|/\Phi$
 and estimate this integral $J$ as follows:
\begin{multline*}
J=
-\int_0^\tau\int_{\mathbb{R}^d}
(\xi u)_\varepsilon u_\varepsilon^k\,  dx\, dt
-2k\int_0^\tau\int_{\mathbb{R}^d}
(a^{ij}\partial_{x_i}\Phi \varrho)_\varepsilon
u_\varepsilon^{k-1}\partial_{x_j}u_\varepsilon\,  dx\, dt
\\
\le
\Bigl(\int_0^\tau \int_{\mathbb{R}^d}
|\xi|^{\beta/2}u
\,  dx\, dt\Bigr)^{2/\beta}
\Bigl(\int_0^\tau  \int_{\mathbb{R}^d}
u_\varepsilon^{k\beta/(\beta -2)+1}\,  dx\, dt
\Bigr)^{(\beta -2)/\beta}
\\
+2k\Bigl(\int_0^\tau\int_{\mathbb{R}^d}
(\eta u)_\varepsilon^2 u_\varepsilon^{k-1}
\, dx\, dt\Bigr)^{1/2}
\Bigl(\int_0^\tau\int_{\mathbb{R}^d}
|\nabla u_\varepsilon|^2 u_\varepsilon^{k-1}
\, dx\, dt\Bigr)^{1/2}.
\end{multline*}
It remains to observe that there hold
the equality
$$
\int_0^\tau \int_{\mathbb{R}^d}
|\xi|^{\beta/2}u \,  dx\, dt
=
\int_0^\tau \int_{\mathbb{R}^d}
|L\Phi|^{\beta/2}\Phi^{1-\beta/2} \varrho
\,  dx\, dt
$$
and the inequality
$$
\int_0^\tau\int_{\mathbb{R}^d}
(\eta u)_\varepsilon^2 u_\varepsilon^{k-1}
\, dx\, dt
\le
 \Bigl(\int_0^\tau\int_{\mathbb{R}^d}
\eta^\beta \Phi \varrho \, dx\, dt \Bigr)^{2/\beta}
\Bigl(\int_0^\tau \int_{\mathbb{R}^d}
u_\varepsilon^{(k\beta)/(\beta-2) +1}\,  dx\, dt
\Bigr)^{(\beta-2)/\beta},
$$
which is verified in the same manner as
in Lemma~\ref{lem3.2}.
Since $s=\beta$, one has $r'\le \beta/2$.
\end{proof}

\begin{example}
{\rm
Suppose that $A$ and $A^{-1}$ are uniformly bounded,
$A$ is uniformly Lipschitzian in~$x$, and for some
$\beta>d+2$, $r>0$,
$\varepsilon>0$, $K>0$ one has
\begin{equation}\label{e3.8}
|b|\in L^\beta(\mu),\
 \exp[(2K+\varepsilon) |x|^r]\in L^{1}(\mu),\
 \sup_{t\in [0,1]}
 \int_{\mathbb{R}^d} \exp(K|x|^r)\varrho(t,x)\, dx<\infty.
\end{equation}
Let
$\sup_{t\in[0,1]} \|b(t,\,\cdot\,)\|_{L^d(\mu_t)}<\infty$.
Finally, let the function $\exp(K|x|^r)\varrho(0,x)$
be bounded and integrable on~$\mathbb{R}^d$.
Then for every $\tau<1$
there is a number $C(\tau)$ such that
$$
\varrho(t,x)\le C(\tau)\exp(-K|x|^r),\quad
(t,x)\in [0,\tau]\times\mathbb{R}^d.
$$
In order to ensure condition (\ref{e3.8}) and the assumptions
on $b$ and $\varrho(0,\,\cdot\,)$ it suffices to have
the estimates
$|b(t,x)|\le C\exp (2K\beta^{-1}|x|^r)$,
$\varrho(0,x)\le C\exp(-K'|x|^r)$ with  $K'>K$  and
\begin{equation}\label{e3.9}
(x,b(t,x)) \le c_1-c_2|x|^r,\
c_2> 2r K\sup_{t,x}\|A(t,x)\|.
\end{equation}
Indeed, let $\Phi\in C^2(\mathbb{R}^d)$,
$\Phi(x)=\exp(K|x|^r)$ при $|x|\ge 1$.
All hypotheses of Theorem \ref{t3.3} are fulfilled.
Under condition (\ref{e3.9}) we pick
 $\delta\in (0,\varepsilon)$ such that one has
the inequality
$
r(2K+\delta)\sup_{t,x}\|A(t,x)\|<c_2,
$
and take a function $V\in C^2(\mathbb{R}^d)$ that
equals $\exp[(2K+\delta)|x|^r]$ if $|x|\ge 1$.
Then, for some~$c$, we have the estimate $LV\le c$.
It follows from \cite{BDPR} that a solution exists
and the norms
$\|V\varrho(t,\,\cdot\,)\|_{L^1(\mathbb{R}^d)}$
are uniformly bounded.
Other assumptions
of Theorem \ref{t3.3} are fulfilled as well.
Similarly, under weaker conditions,
one can obtain a power bound.
}\end{example}

Analogous theorems are valid in the situation of the second theorem
of the previous section.

Most of the work
has been done during visits of the first
and third authors to the
University of Bielefeld.

\end{document}